\documentclass[12pt,reqno,a4paper]{amsart}

\usepackage[T1]{fontenc}
\usepackage{lmodern}
\usepackage{amsmath,amssymb,amsthm,mathtools}
\usepackage[margin=3cm]{geometry}
\usepackage[hidelinks]{hyperref}

\newcommand{\R}{\mathbb R}
\newcommand{\E}{\mathbb E}
\newcommand{\Var}{\operatorname{Var}}
\newcommand{\KL}{\operatorname{H}}
\newcommand{\ind}{\mathbf{1}}

\newtheorem{theorem}{Theorem}[section]
\newtheorem{proposition}[theorem]{Proposition}
\newtheorem{lemma}[theorem]{Lemma}
\newtheorem{corollary}[theorem]{Corollary}
\theoremstyle{remark}
\newtheorem{remark}[theorem]{Remark}

\newcommand{\N}{{\mathbb{N}}}

\numberwithin{equation}{section}

\title[Speed of convergence for propagation of chaos]{The speed of convergence for propagation of chaos in the mean-field Ising model}

\begin{document}
	
\author[Jonas Jalowy]{Jonas Jalowy}
\address[Jonas Jalowy]{Institut f\"ur Mathematik,
	Universit\"at Paderborn,
	Warburgerstr. 100, 33098 Paderborn, Germany}

\email[Jonas Jalowy]{jjalowy@math.uni-paderborn.de}

\author[Matthias L\"owe]{Matthias L\"owe }
\address[Matthias L\"owe]{Fachbereich Mathematik und Informatik,
	Universit\"at M\"unster,
	Einsteinstra\ss e 62,
	48149 M\"unster,
	Germany}
\email[Matthias L\"owe]{maloewe@uni-muenster.de}

	\author[Sara Terveer]{Sara Terveer}
\address[Sara Terveer]{Mathematisches Institut,\\
	Ludwig-Maximilians-Universit\"at M\"unchen,\\
	Theresienstraße 39,
	80333 M\"unchen,
	Germany}
\email[Sara Terveer]{terveer@math.lmu.de}

\subjclass[2020]{Primary: 82B05, Secondary: 82B20, 82B27,60B10}
\keywords{Curie–Weiss model, propagation of chaos, relative entropy, chi-square divergence, pair correlations, critical scaling.}

\begin{abstract}
We derive sharp relative-entropy asymptotics for increasing propagation of chaos in the Curie--Weiss model at zero external field. 
For every fixed inverse temperature $0<\beta<1$, the relative entropy of
a $k$-spin marginal with respect to the symmetric product measure is
$$
\frac{\beta^2}{4(1-\beta)^2}\frac{k(k-1)}{N^2}+o\left(\frac{k^2}{N^2}\right),\qquad k=o(N).
$$
At the critical inverse temperature $\beta=1$, the rate changes to
$\alpha_{\mathrm c}^2 k(k-1)/(4N)+o(k^2/N)$ in the regime $k=o(\sqrt N)$,
where $\alpha_{\mathrm c}$ is the second moment of the limiting quartic
law of the rescaled magnetization. We also obtain asymptotics for sequences $\beta_N\uparrow1$ that interpolate between these two regimes. 

Crucially, the proof reveals that the pair correlations determine  the leading correction to the likelihood ratio, providing precise control of both the chi-square divergence and relative entropy.
\end{abstract}

\maketitle

\section{Introduction and main results}

Propagation of chaos is one of the fundamental paradigms of mean-field
theory; see, for instance, \cite{propChaos,Kac_foundations,sznitman_poc}. For a sequence of symmetric
probability measures on an $N$-particle space, it means that, for every
fixed $k$, the $k$-particle marginal converges to the $k$-fold product of a
limiting one-particle measure. If several macroscopic phases coexist,
the limiting object may instead be a mixture of product measures.

Ben Arous and Zeitouni \cite{BAZ_chaos} introduced increasing propagation
of chaos, in which the number $k=k(N)$ of observed particles may diverge
with the system size. They prove propagation of chaos for marginals of size $k=o(N)$, when
the associated free-energy functional has a unique non-degenerate
minimizer, discuss several minimizers and degenerate minima, and treat the
high-temperature Curie--Weiss model as a fundamental example. See also \cite{propChaos,Papa,Wilms,JKLM,JKL25,GLS26} for related results.

Ben Arous and Zeitouni \cite{BAZ_chaos} first bound the relative entropy of the full Gibbs measure with respect to the  product measure. In the non-degenerate case, this entropy remains bounded as $N\to  \infty$. Csisz\'ar's entropy inequality  \cite{Csiszar67}
then yields an upper bound of order $k/N$ for the relative entropy
of the $k$-particle marginal. This bound need not be sharp. For a class of continuous
mean-field Gibbs measures on $(\mathbb R^d)^N$, Lacker \cite{lacker2021quantitative,Lacker}
obtained the sharper rate $\mathcal O (k^2/N^2)$ in a high-temperature regime. This automatically raises the question how fast the speed of convergence in propagation of chaos is in discrete systems, where the methods of \cite{Lacker} appear to be inapplicable.  Here, we answer this question for a particular discrete model, the Curie-Weiss model. More precisely, we determine the sharp quadratic rate $\mathcal O (k^2/N^2)$ for the speed of convergence, throughout the entire regime $k=o(N)$. We also identify the change of rate at criticality and the transition when the inverse temperature approaches its critical value.

\subsection{The model}
Let us start by describing the model we have in mind: 
For $N\in\mathbb N$ and $\beta>0$, the Curie--Weiss measure without
external field is given by the Gibbs measure. 
\begin{align}\label{eq:mu_N}
\mu_{N,\beta}(\sigma):=\frac{1}{Z_{N,\beta}}
\exp\left(\frac{\beta}{2N}\left(\sum_{i=1}^N\sigma_i\right)^2\right),
\qquad \text{where} \quad \sigma\in\{-1,+1\}^N, 
\end{align}
where $Z_{N,\beta}$ is the normalizing constant. The marginal of the reference measure is given by the symmetric Bernoulli measure
\[
\nu_0:=\frac12(\delta_{-1}+\delta_{+1}).
\]
Further, we write
$\bar\mu_{N,\beta}^{(k)}$ for the marginal of $\mu_{N,\beta}$ on any $k$ coordinates. The specific choice of coordinates is irrelevant by exchangeability. Obviously, $\bar\mu_{N,\beta}^{(1)}$ is the same as $\nu_0$ for every $N$.
The results in \cite{BAZ_chaos,JKLM} show that
$\bar\mu_{N,\beta}^{(k)}$ approaches $\nu_0^{\otimes k}$
for fixed $0<\beta<1$ when $k=o(N)$, and for $\beta=1$
when $k=o(\sqrt N)$.
Our aim is to determine the leading asymptotic value
of the relative entropy in these regimes.  

This already shows that the model exhibits a phase transition at the critical inverse temperature $\beta_c=1$, which is detectable 
via the spin sum $\tau_N:=\sum_{i=1}^N\sigma_i$ and the magnetization $\tau_N/N$. For $\beta\le 1$, the magnetization converges to zero in $\mu_{N,\beta}$ probability, whereas for $\beta>1$, it concentrates around two nonzero values of opposite sign.
We refer to \cite{EllisEntropyLargeDeviationsAndStatisticalMechanics,
BovierSMoDS,Velenik_book} for further background on the model. The fluctuations of the magnetization
explain why two different propagation-of-chaos scales arise:
For every fixed $\beta<1$,
\begin{align}\label{eq:HT-CLT}
\frac 1 {\sqrt N} \tau_N=\frac 1 {\sqrt N} \sum_{i=1}^N \sigma_i \Rightarrow \mathcal N\left(0,\frac1{1-\beta}\right)
\end{align}
where $\mathcal N(\mu, \sigma^2)$ denotes the normal distribution with 
expectation $\mu$ and variance $\sigma^2$.
At $\beta=1$, in accordance with the CLT at $\beta <1$ the fluctuations are larger:
$N^{-3/4}\tau_N$ converges in distribution to a random variable $Y$ with density $\propto e^{-y^4/12}$. 
See also \cite{EL10,Chatterjee_Shao} for the critical fluctuation result.
The second moment of this law will determine the coefficient in our
critical entropy asymptotics:
\begin{equation}\label{eq:alphac}
\alpha_{\mathrm c}:=\E[Y^2]
=\frac{\displaystyle\int_{\mathbb R}y^2e^{-y^4/12} dy}
{\displaystyle\int_{\mathbb R}e^{-y^4/12} dy}
=2\sqrt3 \frac{\Gamma(3/4)}{\Gamma(1/4)},
\end{equation}
where the last equality follows
by the substitution $t=y^4/12$ on $(0,\infty)$.

For probability measures $P$ and $Q$ on the same measurable space,
with $P\ll Q$, we quantify the deviation of $P$ from $Q$
by the relative entropy 
\begin{align}
\KL(P\mid Q):=\int\log\left(\frac{dP}{dQ}\right) dP
\end{align}
and the chi-square divergence
\begin{align}
\chi^2(P,Q):=\int\left(\frac{dP}{dQ}-1\right)^2 dQ.
\end{align}
For completeness, we set $\KL(P\mid Q)=+\infty$ if $P\not\ll Q$. We have that $H\le \chi^2$, which follows immediately from $x\log x\le x(x-1)$.

For the marginals considered below, we will prove in Lemma \ref{lem:entropy-quadratic} that relative
entropy is asymptotically half the chi-square divergence.
Our proofs therefore first determine the chi-square divergence
and then deduce the relative entropy asymptotics.

\subsection{Main results}

The first main result gives the leading asymptotic of the propagation of chaos in the
high-temperature regime.

\begin{theorem}\label{thm:entropy-high-temperature}
Let $0<\beta<1$, and let $k=k(N)$ satisfy $1\le k\le N$ and $k=o(N)$.
Then
\[
\chi^2\left(\bar\mu_{N,\beta}^{(k)},\nu_0^{\otimes k}\right)
=\frac{\beta^2}{2(1-\beta)^2}\frac{k(k-1)}{N^2}
+o\left(\frac{k^2}{N^2}\right),
\]
and
\[
\KL\left(\bar\mu_{N,\beta}^{(k)}\mid\nu_0^{\otimes k}\right)
=\frac12\chi^2\left(\bar\mu_{N,\beta}^{(k)},\nu_0^{\otimes k}\right)
+o\left(\frac{k^2}{N^2}\right)=\frac{\beta^2}{4(1-\beta)^2}\frac{k(k-1)}{N^2}+o\left(\frac{k^2}{N^2}\right).
\]
\end{theorem}
 The factor $\binom k2=k(k-1)/2$ is the number of distinct spin-pairs among the first $k$ spins; in particular, the error vanishes for $k=1$ as expected. 
 
 \begin{remark}[Related literature]
 	As mentioned above, Lacker \cite[Corollary 2.7]{lacker2021quantitative,Lacker} obtained related $ \mathcal O (k^2/N^2)$ upper bounds for the relative entropy of $k$-particle marginals of certain continuous high-temperature mean-field Gibbs measures. His proof relies on a differentiation identity  \cite[Lemma 3.1]{lacker2021quantitative} for the logarithmic derivative of the density, which is inapplicable in discrete systems, and does not optimize his bound \cite[(2.10)]{lacker2021quantitative}. In contrast, the leading coefficient in Theorem \ref{thm:entropy-high-temperature} is explicit.

 Let us also compare our results to \cite[Equation (27)]{Wilms}, who consider discrete mean-field spin systems and the regime of fixed block fractions $ k/N\sim t\in(0,1)$. In our notation, their saddle-point approximation gives
\begin{align}\label{eq:Wilms} \KL\left(\bar\mu_{N,\beta}^{(tN)}\mid\nu_0^{\otimes tN}\right)=\frac12\Bigl(\frac{\beta t}{1-\beta}-\log\big(1+\frac{\beta t}{1-\beta}\big)\Bigr)+\mathcal O (N^{-1}),\qquad N\to\infty.\end{align}
A formal, non-rigorous expansion at $t=k/N\to 0$ recovers the constants of Theorem 1.1. 

We establish rigorous asymptotics throughout the joint regime $k=o(N)$, including bounded $k$ and, show that even a uniform $\mathcal O (N^{-1})$ remainder in \eqref{eq:Wilms} would be too large to capture the leading term in the regime $k=o(\sqrt N)$. We also determine how it changes at and near criticality, where the pair correlations are larger, leading to divergent entropies in \cite{Wilms} and to even smaller admissible block sizes for our results.
\end{remark}

\begin{theorem}\label{thm:critical-entropy}
Let $\beta=1$, and let $k=k(N)$ satisfy $1\le k\le N$ and $k=o(\sqrt N)$.
Then
\[
\chi^2\left(\bar\mu_{N,1}^{(k)},\nu_0^{\otimes k}\right)
=\frac{\alpha_{\mathrm c}^2}{2}\frac{k(k-1)}{N} +o\left(\frac{k^2}{N}\right),
\]
and
\[
\KL\left(\bar\mu_{N,1}^{(k)}\mid\nu_0^{\otimes k}\right)
=\frac{\alpha_{\mathrm c}^2}{4}\frac{k(k-1)}{N} +o\left(\frac{k^2}{N}\right).
\]
\end{theorem}

The windows $k=o(N)$ and $k=o(\sqrt N)$ imply qualitative propagation
of chaos and they are optimal, see \cite{Papa,Wilms,JKLM,JKL25,GLS26}.

To describe the transition between these two results, we will consider an N-dependent temperature 
$\beta_N\uparrow1$. The high-temperature central limit theorem \eqref{eq:HT-CLT} suggests the
magnetization scale $[N(1-\beta_N)]^{-1/2}$, whereas the critical
scale is $N^{-1/4}$. These scales agree when $1-\beta_N$ is of
order $N^{-1/2}$. We therefore introduce the parameter
\begin{equation}\label{eq:def_lambda}
\lambda_N:=\sqrt N(1-\beta_N)
\end{equation}
to describe the transition between the two fluctuation regimes.
When $\lambda_N$ stays bounded (but does not converge to $0$), both the quadratic and quartic terms
contribute. When $\lambda_N\to  \infty$, the quadratic term determines
the fluctuation scale. This leads to the following interpolating
coefficient.

\begin{theorem}\label{thm:near-critical}
For $\lambda\ge0$, define
\[
\alpha(\lambda):=
\frac{\displaystyle\int_{\mathbb R}y^2
\exp\left(-\frac\lambda2y^2-\frac1{12}y^4\right) dy}
{\displaystyle\int_{\mathbb R}
\exp\left(-\frac\lambda2y^2-\frac1{12}y^4\right) dy}.
\]
Then $\alpha:[0,\infty)\to  (0,\alpha_{\mathrm c}]$ is continuous and
strictly decreasing, with
\[
\alpha(0)=\alpha_{\mathrm c} \qquad 
\text{and} \qquad 
\alpha(\lambda)\sim\frac1\lambda
\quad\text{as }\lambda\to  \infty.
\]
Let $\beta_N\uparrow1$, define $\lambda_N$ as in \eqref{eq:def_lambda}, and let $1\le k=k(N)\le N$ satisfy
\[
\frac{k\alpha(\lambda_N)}{\sqrt N}\longrightarrow  0.
\]
Then
\begin{equation}\label{eq:chi2window}
\chi^2\left(\bar\mu_{N,\beta_N}^{(k)},\nu_0^{\otimes k}\right)
=\frac{\alpha(\lambda_N)^2}{2}\frac{k(k-1)}{N}
+o\left(\frac{k^2\alpha(\lambda_N)^2}{N}\right),
\end{equation}
and
\begin{equation}\label{eq:KLwindow}
\KL\left(\bar\mu_{N,\beta_N}^{(k)}\mid\nu_0^{\otimes k}\right)
=\frac{\alpha(\lambda_N)^2}{4}\frac{k(k-1)}{N}
+o\left(\frac{k^2\alpha(\lambda_N)^2}{N}\right).
\end{equation}
\end{theorem}

For bounded $\lambda_N$, the critical condition $k=o(\sqrt N)$ remains the relevant one. If $\lambda_N\to  \infty$, the condition becomes $k=o(N(1-\beta_N))$, and the leading entropy is
$k(k-1)/(4N^2(1-\beta_N)^2)$. Since $\beta_N\to  1$, this agrees with the coefficient suggested by the fixed-temperature result. Among the many works on finite-size scaling near the critical temperature, we mention the magnetization scaling form \cite{Botet}, the critical window description via the Ginzburg parameter \cite[Eqn (17)]{Colonna}, and fluctuation results  \cite{EL10}.

As a byproduct of the pair-correlation expansion in our proof, we also obtain a first-order correction for
smooth observables $g\in C_b^2(\mathbb R)$ of the block magnetization $\tau_k:=\sum_{i=1}^k\sigma_i$. For fixed $\beta<1$ and $k\to\infty$ with $k=o(N)$, we will prove that
\[
\frac Nk\left(
\mathbb E_{\bar\mu_{N,\beta}^{(k)}}
\left[g\left(\frac{\tau_k}{\sqrt{k}}\right)\right]
-
\mathbb E_{\nu_0^{\otimes k}}
\left[g\left(\frac{\tau_k}{\sqrt{k}}\right)\right]
\right)
\longrightarrow
\frac{\beta}{2(1-\beta)}\mathbb E[g''(Z)],
\]
where $Z\sim\mathcal N(0,1)$, see Corollary \ref{cor:cor!} below.

\subsection{Proof ideas}
Each individual spin has distribution $\nu_0$ for every $N$.
The discrepancy between the $k$-spin marginal and
$\nu_0^{\otimes k}$ is therefore entirely due to dependence
between the spins. It turns out that at first order, only pair interactions matter. To study this dependence in terms of the relative entropy or $\chi^2$, we 
need to control the likelihood ratio (or, Radon-Nikodym density) 
\[
L_{N,k}:=
\frac{d\bar\mu_{N,\beta}^{(k)}}{d\nu_0^{\otimes k}}.
\] 
 Let us describe the common argument for all our main results, which we derive rigorously in Section \ref{sec:common}. We begin with the standard Hubbard--Stratonovich transformation to express
the Gibbs measure as a mixture of product measures similar to the findings of \cite{deb_mukherjee,JKLM,JKL25}. More precisely, Proposition \ref{prop:mixture} reveals some mixing
variable $X\sim\gamma_{N,\beta}$ such that 
\[
L_{N,k}(\sigma)=\E_{\gamma_{N,\beta}}\Big[
\prod_{i=1}^k \Big(1+\sigma_i\tanh\Big(\sqrt{\frac{\beta}{N}} X\Big)\Big)
\Big].
\]
Thus, conditionally on $X$, the spins $\sigma_i$ are independent.
Expanding the product yields an expansion in products of spins, whose coefficients are the spin correlations
\begin{align}\label{eq:qNr}
q_{N,r}:=\E_{\mu_{N,\beta}}[\sigma_1\cdots\sigma_r]=\E_{\gamma_{N,\beta}}\big[\tanh^r\big(\sqrt{\tfrac{\beta}{N}}X\big)\big],
\end{align}
see~\eqref{eq:expectation_sigmas}. All odd terms vanish by spin-flip symmetry, hence we shall find
\begin{align}\label{eq:expansion_overview}
L_{N,k}(\sigma)
=1+q_{N,2}\sum_{1\le i<j\le k}\sigma_i\sigma_j
+\text{terms involving four or more spins},
\end{align}
leading to the common argument that is summarized in Proposition~\ref{prop:pair-dominance}: Assuming uniform bounds on the higher correlations $q_{N,2r}$, the terms involving four or more spins are negligible whenever $kq_{N,2}\to0$.  On the other hand, all involved spin products $\sigma_A=\prod_{i\in A}\sigma_i$ are orthogonal under the product measure $\nu_0^{\otimes k}$, since the spins are independent and centered. Averaging our expansion, we see that only the pair-correlations contribute to the chi-square divergence with
\[
\chi^2\left(\bar\mu_{N,\beta}^{(k)},\nu_0^{\otimes k}\right)
\sim
\binom{k}{2}q_{N,2}^2,\qquad k\ge 2.
\]  
The connection to the entropy via $H\sim\tfrac 1 2 \chi^2$ follows from a Taylor approximation of $x\log x$, see Lemma \ref{lem:entropy-quadratic} below. Therein, we will make use of the Bonami inequality \cite[Theorem~9.21]{ODonnell14}, which bounds higher moments of homogeneous polynomials in the spins by lower ones (Boolean hypercontractivity), see also \cite[Eqn (17)]{sessak2009small} for a related result.

It then remains to solely  study pair correlations $q_{N,2}$ and the assumed uniform correlation bound in each of the three temperature regimes separately in Sections~3--5.
 The distinction between them is contained in the mixing law: For fixed $\beta<1$, it is approximatively  Gaussian and we get $q_{N,2}\sim\beta/((1-\beta)N)$. At $\beta=1$, its fluctuation scale
 changes, as usual, and its limit has a quartic exponential density, yielding
 $q_{N,2}\sim\alpha_{\mathrm c}/\sqrt N$. When $\beta_N\uparrow1$,
 the limiting quadratic and quartic terms give  $q_{N,2}\sim\alpha(\lambda_N)/\sqrt N$.

\section*{Acknowledgements and AI-statement}
JJ has been financially supported by the DFG priority program SPP 2265 \emph{Random Geometric Systems}.
ML's research was funded by the Deutsche Forschungsgemeinschaft (DFG, German Research Foundation) under Germany's Excellence Strategy EXC 2044/2 - 390685587, Mathematics M\"unster: \emph{Dynamics-Geometry-Structure}.
ST's research was funded by the Deutsche Forschungsgemeinschaft (DFG, German Research Foundation)  -- TRR 352 ``Mathematics of Many-Body Quantum Systems and Their Collective Phenomena''
-- Project-ID 470903074.

The authors obtained the main results and identified the leading contribution $\frac{q_{N,2}}2(\tau_k^2-k)$ before AI was available. Later on, these findings have been prompted to ChatGPT~6, which (after extensive prompting), assisted in developing the less technical proof presented here. AI helped with step-by-step drafting and proof-checks, and all AI-assisted material was subsequently carefully reviewed, verified, and humanized.

\section{The common core: reduction to pair correlations}\label{sec:common}


As announced above, we start with the usual Hubbard-Startonovich transformation. It replaces the interaction of all spins by the interactions with a (real-valued) field $X$ having distribution 
\begin{align}\label{eq:gamma}
\gamma_{N,\beta}(dx):=
\frac{e^{F_{N,\beta}(x)}}{\displaystyle\int_{\mathbb R} e^{F_{N,\beta}(y)} d y}dx,
\end{align}
for \begin{align}\label{eq:F_N}
F_{N,\beta}(x):=-\frac{x^2}{2}+N\log\cosh\Big(\sqrt{\frac\beta N}x\Big).
\end{align}
Here, we keep the dependence on $\beta$ explicit, since  $\beta$ may depend on $N$ later on, but in order not to overload notation, we suppress the dependence for $L_{N,k}$ and the
correlations $q_{N,r}$.
\begin{proposition}\label{prop:mixture}
Let $X\sim\gamma_{N,\beta}$. For every $1\le k\le N$,
\begin{equation}\label{eq:representation_L}
L_{N,k}(\sigma) =\E_{\gamma_{N,\beta}}\left[
\prod_{i=1}^k\bigl(1+\sigma_i\tanh\left(\sqrt{\frac\beta N} X\right)\bigr)\right]
\qquad \text{for } \quad \sigma\in\{-1,+1\}^k.
\end{equation}
Equivalently, $\bar\mu_{N,\beta}^{(k)}$ is the mixture of product
measures whose one-spin probabilities, conditional on $X=x$, are
$\frac12\big(
1+\sigma_i\tanh\big(\sqrt{\frac{\beta}{N}}x\big)
\big)$ with mean $\tanh\big(\sqrt{\frac{\beta}{N}}x\big)$.
\end{proposition}

\begin{proof}
The Hubbard--Stratonovich (or just completion of the square) identity reads
\[
\exp\left(\frac{\beta s^2}{2N}\right) =\frac1{\sqrt{2\pi}}\int_{\mathbb R}
\exp\left(-\frac{x^2}{2}+\sqrt{\frac\beta N} xs\right) dx.
\] 
Applying it to $s=\sum_{i=1}^N\sigma_i$, using the definition \eqref{eq:mu_N} of the Gibbs measure and summing over the spins gives
\[
\begin{aligned}
	Z_{N,\beta}
	&=\frac1{\sqrt{2\pi}}\int_{\mathbb R}e^{-x^2/2}
	\sum_{\sigma\in\{-1,+1\}^N}
	\exp\left(\sqrt{\frac{\beta}{N}} x\sum_{i=1}^N\sigma_i\right) dx\\
	&=\frac{2^N}{\sqrt{2\pi}}\int_{\mathbb R}
	e^{-x^2/2}\cosh^N\left(\sqrt{\frac{\beta}{N}} x\right) dx =\frac{2^N}{\sqrt{2\pi}}\int_{\mathbb R}e^{F_{N,\beta}(x)} dx.
\end{aligned}
\] 
For a fixed $k$-spin configuration, summation over the other $N-k$
spins therefore yields
\[
\bar\mu_{N,\beta}^{(k)}(\sigma)
=2^{-k}\int_{\mathbb R}
\prod_{i=1}^k \frac{\exp(\sqrt{\beta/N} x\sigma_i)} {\cosh(\sqrt{\beta/N} x)} \gamma_{N,\beta}(dx).
\]
For $\sigma_i\in\{-1,+1\}$, we have
$e^{t\sigma_i}/\cosh t=1+\sigma_i\tanh t$.
Dividing by $\nu_0^{\otimes k}(\sigma)=2^{-k}$ proves the formula.
\end{proof}
Higher moments of $\tanh(\sqrt{\frac{\beta}{N}}X)$ have the probabilistic meaning of describing the spin
correlations under the Gibbs measure. 
For $A\subseteq[k]:=\{1,\ldots,k\}$, write
$\sigma_A:=\prod_{i\in A}\sigma_i$, with $\sigma_\emptyset:=1$. By conditional independence from Proposition \ref{prop:mixture} we have for any $1\le r\le N$,
\begin{align}\label{eq:expectation_sigmas}
q_{N,r}=	\E_{\mu_{N,\beta}}[\sigma_1\cdots\sigma_r]&=
	\E_{\gamma_{N,\beta}}\Big[\E[\sigma_1\cdots\sigma_r\mid X]\Big]\\
	&=\E_{\gamma_{N,\beta}}\Big[\prod_{i=1}^r\E[\sigma_i\mid X]\Big] =\E_{\gamma_{N,\beta}}\Big[\tanh^r\Big(\sqrt{\frac{\beta}{N}}X\Big)\Big]=:q_{N,r},\nonumber
\end{align}
as claimed in \eqref{eq:qNr}. In particular, $r=2$ gives the pair-correlation from \eqref{eq:expansion_overview}. The right hand side of \eqref{eq:expectation_sigmas} may even serve as a definition for $q_{N,r}$ in the case of $r>\lfloor N/2\rfloor$ and $r=0$.

\begin{remark}\label{rem:naiv_bound}
Proposition~\ref{prop:mixture} suggests to naively bound the relative
entropy by averaging the entropies of the conditional product measures.
By convexity of relative entropy, its additivity for
product measures, and the inequality $\KL\le\chi^2$,
we obtain the following valid upper bound 
\begin{align*}
	\KL\left(\bar\mu_{N,\beta}^{(k)}
	\mid\nu_0^{\otimes k}\right)
	&\le
	\E_{\gamma_{N,\beta}}\left[\KL\left(\bar\mu_{N,\beta}^{(1)}(\cdot|X)^{\otimes k}\mid\nu_0^{\otimes k}\right)
	\right]=k\E_{\gamma_{N,\beta}}\left[\KL(\bar\mu_{N,\beta}^{(1)}(\cdot|X)\mid\nu_0)\right]\\
	&\le k\E_{\gamma_{N,\beta}}\big[\tanh^2\big(\sqrt{\tfrac{\beta}{N}}X\big)\big]
	=kq_{N,2}.
\end{align*} 
However, this bound is not sharp as it does not capture the cancellation of the one-spin
biases under the symmetric mixing law.
For a fixed value of the common field, the spins are independent
but generally biased with $\KL(\bar\mu_{N,\beta}^{(1)}(\cdot\mid X)\mid\nu_0)>0$. Mixing over the symmetric law of $X$ removes these one-spin biases, giving $\bar\mu_{N,\beta}^{(1)}=\nu_0$. Already for $k=1$, the true relative entropy is exactly zero, whereas the above bound is strictly positive.
\end{remark}

To obtain a sharper estimate, we therefore directly expand the averaged likelihood ratio \eqref{eq:representation_L}
 
\begin{proposition}\label{prop:walsh}
For every $k\le N$ we have
\begin{equation}\label{eq:LNWalsh}
L_{N,k}(\sigma)-1
=\sum_{r=1}^{\lfloor k/2\rfloor}q_{N,2r}
\sum_{\substack{A\subseteq[k]\\|A|=2r}}\sigma_A.
\end{equation}
Moreover,
\begin{equation}\label{eq:chi-square-exact}
\chi^2\left(\bar\mu_{N,\beta}^{(k)},\nu_0^{\otimes k}\right)
=\sum_{r=1}^{\lfloor k/2\rfloor}\binom{k}{2r}q_{N,2r}^2.
\end{equation}
\end{proposition}

\begin{proof}
Expanding the product in \eqref{eq:representation_L} gives
\[
L_{N,k}(\sigma)=\sum_{A\subseteq[k]}q_{N,|A|}\sigma_A.
\]
The empty set contributes one. Correlation terms $q_{N,|A|}$ with odd $|A|$ vanish, since the mixing law $\gamma_{N,\beta}$ is symmetric and $\tanh (\sqrt{\frac{\beta}{N}}x )$ is odd.  
This proves \eqref{eq:LNWalsh}.

Under $\nu_0^{\otimes k}$, the spins are independent and centered, implying the claimed orthogonality:
If $A=B$, then $\E[\sigma_A\sigma_B]=1$ and if $A\ne B$, then some spin occurs exactly once in the product and thus $\E[\sigma_A\sigma_B]=0$.
Squaring \eqref{eq:LNWalsh} and taking expectations therefore gives
\eqref{eq:chi-square-exact}.
\end{proof}

The expansion in \eqref{eq:LNWalsh} is sometimes also called Fourier--Walsh expansion (of the likelihood ratio),  see \cite{ODonnell14}. The term indexed by
$r$ involves products of $2r$ spins. The first term $r=1$ is a sum over
pairs and can also be written in terms of the block sum $\tau_k=\sum_{i=1}^k\sigma_i$; 
\begin{equation}\label{eq:pair-block-sum}
q_{N,2}\sum_{1\le i<j\le k}\sigma_i\sigma_j
=\frac{q_{N,2}}2(\tau_k^2-k).
\end{equation}
This term is centered under the product distribution, but not under the Gibbs measure. Before we show that it dominates higher correlation terms, 
 we establish a criterion for passing from chi-square
asymptotics to relative entropy asymptotics.

\begin{lemma}\label{lem:entropy-quadratic}
Let $P_N\ll Q_N$ be probability measures, and write $L_N=dP_N/dQ_N$.
If
\[
\chi^2(P_N,Q_N):=\E_{Q_N}[(L_N-1)^2]\to   0
\qquad \text{and} \qquad 
\E_{Q_N}[(L_N-1)^4] =o\left(\E_{Q_N}[(L_N-1)^2]\right),
\]
then
\[
\KL(P_N\mid Q_N) =\frac12\chi^2(P_N,Q_N) +o\left(\E_{Q_N}[(L_N-1)^2]\right).
\]
\end{lemma}

\begin{proof}
Since $L_N$ is a probability density with respect to $Q_N$, we have 
\[
\KL(P_N\mid Q_N)
=\E_{Q_N}\left[L_N\log L_N-(L_N-1)\right].
\]
We want to compare the integrand with
$\frac12(L_N-1)^2$. Writing $L_N-1= u$, the difference is
\[
(1+u)\log(1+u)-u-\tfrac12u^2.
\]
By Taylor's formula at $u=0$ we obtain for all
$\varepsilon>0$, that there is a $\delta\in(0,1)$ such that
\begin{align}
|(1+u)\log(1+u)-u-\tfrac12u^2|
\le\varepsilon u^2
\qquad\text{if }|u|\le\delta.\label{eq:Taylor}
\end{align}
On the other hand, a short analysis shows that for all $u\ge-1$
\[
0\le(1+u)\log(1+u)-u\le u^2,
\]
where we set $0 \log 0=0$.
Hence, $|(1+u)\log(1+u)-u-\frac12u^2|\le u^2/2$, even for $|u|>\delta$. We split according to whether $|L_N-1|\le\delta$ and use
\[
(L_N-1)^2\ind_{\{|L_N-1|>\delta\}} \le\delta^{-2}(L_N-1)^4,
\]
to obtain
\[
\left|\KL(P_N\mid Q_N)-\frac12\E_{Q_N}[(L_N-1)^2]\right|
\le\varepsilon\E_{Q_N}[(L_N-1)^2] +\frac1{2\delta^2}\E_{Q_N}[(L_N-1)^4].
\]
Now fix $\varepsilon>0$ and choose $\delta$ as in \eqref{eq:Taylor}.
Whenever the second moment is positive, we divide by it to obtain
\[
\frac{\left|\KL(P_N\mid Q_N)-\frac12\E_{Q_N}[(L_N-1)^2]\right|}
{\E_{Q_N}[(L_N-1)^2]}
\le\varepsilon+\frac1{2\delta^2}\frac{\E_{Q_N}[(L_N-1)^4]}{\E_{Q_N}[(L_N-1)^2]}.
\]
By assumption, the second term on the right hand side tends to zero. This holds for every $\varepsilon>0$, so the relative error tends to zero, proving the claimed asymptotic expansion.
If the second moment is zero, then $L_N=1$ $Q_N$-almost surely
and hence $P_N=Q_N$. In this case, both the relative entropy
and the second moment vanish.
\end{proof}

The next proposition provides the pair-correlation criterion for the common estimates needed in all three
regimes in the theorems.
\begin{proposition}\label{prop:pair-dominance}
Let $\beta=\beta_N$ depend on $N$, if desired, and suppose
$2\le k=k(N)\le N$. Assume that $q_{N,2}>0$, that $kq_{N,2}\to   0$,
and that there is a constant $C$ (independent of $N$ and $r$) such that
\begin{equation}\label{eq:common-moment-assumption}
q_{N,2r}\le(Cr q_{N,2})^r,
\qquad 1\le r\le\lfloor k/2\rfloor,
\end{equation}
for all sufficiently large $N$. Then
\begin{equation}\label{eq:common-chi-asymptotics}
\chi^2\left(\bar\mu_{N,\beta}^{(k)},\nu_0^{\otimes k}\right)
=\binom{k}{2}q_{N,2}^2+O\left((kq_{N,2})^4\right),
\end{equation}
and
\begin{equation}\label{eq:common-entropy-asymptotics}
\KL\left(\bar\mu_{N,\beta}^{(k)}\mid\nu_0^{\otimes k}\right)
=\frac12\binom{k}{2}q_{N,2}^2
+o\left(k^2q_{N,2}^2\right).
\end{equation}
Furthermore, 
\begin{equation}\label{eq:common-L2-remainder}
\E _{\nu_0^{\otimes k}} \left[\Big(L_{N,k}-1-\frac{q_{N,2}}2(\tau_k^2-k)\Big)^2\right]
=O\left((kq_{N,2})^4\right),
\end{equation}
and
\begin{equation}\label{eq:common-fourth-moment}
\E _{\nu_0^{\otimes k}} \left[ \left(L_{N,k}-1\right)^4\right]=O\left((kq_{N,2})^4\right) .
\end{equation}
\end{proposition}

The assumptions have two distinct roles. The condition
$kq_{N,2}\to   0$ makes the total pair contribution small. 
The bound on $q_{N,2r}$ allows us to control the sum over all
terms $r\ge2$ in the expansion \eqref{eq:LNWalsh}, each contributing $\binom{k}{2r}q_{N,2r}^2$ to the chi-square divergence.
 Together with $kq_{N,2}\to   0$, it makes this sum
negligible compared with the pair contribution
$\binom{k}{2}q_{N,2}^2$.

\begin{proof}
We first isolate the pair term in \eqref{eq:chi-square-exact} corresponding to $r=1$.
For $r\ge2$, \eqref{eq:common-moment-assumption} and
$\binom{k}{2r}\le(ek/(2r))^{2r}$ imply
\[
\binom{k}{2r}q_{N,2r}^2
\le\left(\frac{ek}{2r}\right)^{2r}(Crq_{N,2})^{2r} \le (Ckq_{N,2})^{2r},
\]
where constants may change from line to line. Since $kq_{N,2}\to  0$,
summing this geometric bound gives
\begin{align}\label{eq:geo_sum}
\sum_{r=2}^{\lfloor k/2\rfloor}\binom{k}{2r}q_{N,2r}^2\le (Ckq_{N,2})^{4} \sum_{r=0}^{\lfloor k/2\rfloor-2} (Ckq_{N,2})^{2r}
=O\left((kq_{N,2})^4\right).
\end{align}
This proves \eqref{eq:common-chi-asymptotics} by Proposition \ref{prop:walsh}. 

On the other hand, if we remove the pair component \eqref{eq:pair-block-sum} from \eqref{eq:LNWalsh}
and take the second moment of the rest, then orthogonality leads to the same sum \eqref{eq:geo_sum}. This shows \eqref{eq:common-L2-remainder}.

In order to estimate the fourth moment, we use the Bonami inequality (Boolean hypercontractivity):
For a Walsh polynomial of degree at most $d$, its $L^4$-norm is at
most $3^{d/2}$ times its $L^2$-norm; see
\cite[Theorem~9.21]{ODonnell14}. We apply this to each degree $d=2r$ in \eqref{eq:LNWalsh}
separately. 
Orthogonality of the Walsh polynomials and assumption \eqref{eq:common-moment-assumption} yield
\[
\E _{\nu_0^{\otimes k}} \Big[ \Big( q_{N,2r}\sum_{\substack{A\subseteq[k]\\|A|=2r}}
\sigma_A\Big)^4\Big]^{1/4} \le3^r\binom{k}{2r}^{1/2}q_{N,2r}\le 
3^r\left(\frac{ek}{2r}\right)^r(Crq_{N,2})^r \le(Ckq_{N,2})^r.
\]
The sum of these terms in the claim \eqref{eq:common-fourth-moment} is bounded by Minkowski's inequality by summing over all these $L^4(\nu_0^{\otimes k})$ norms and taking the fourth root, for which another geometric sum finishes the proof.

Finally, $k\ge2$ implies $\binom{k}{2}\asymp k^2$. The chi-square estimate \eqref{eq:common-chi-asymptotics} thus gives
\[
\E_{\nu_0^{\otimes k}}[(L_{N,k}-1)^2] \asymp(kq_{N,2})^2\to   0.
\]
The fourth moment is $\mathcal O ((kq_{N,2})^4)$ and is therefore negligible
compared with the second moment. Lemma~\ref{lem:entropy-quadratic}
proves \eqref{eq:common-entropy-asymptotics}.
\end{proof}

The remaining task in each regime is now to calculate $q_{N,2}$ and verify \eqref{eq:common-moment-assumption}. Before we do so, let us record the following corollary which may be of independent interest.
\begin{corollary}\label{cor:testfunction}
	Under the assumptions of Proposition~\ref{prop:pair-dominance}, let
	$f_{N,k}:\{-1,+1\}^k\to\mathbb R$ satisfy
	$\sup_N\mathbb E_{\nu_0^{\otimes k}}[f_{N,k}^2]<\infty$.
	Then,
	\begin{align}
	\frac{
		\mathbb E_{\bar\mu_{N,\beta}^{(k)}}[f_{N,k}]
		-\mathbb E_{\nu_0^{\otimes k}}[f_{N,k}]
	}{kq_{N,2}}
	=
	\frac{1}{2k}
	\mathbb E_{\nu_0^{\otimes k}}
	\bigl[(\tau_k^2-k)\cdot f_{N,k}\bigr]+o(1).
	\end{align} 	
	\end{corollary}
	Thus, to first order, the deviation of any $L^2$-observable from its product-measure
	value is determined by its correlation with the squared block magnetization.
	
	\begin{proof}
It follows from Cauchy--Schwarz and \eqref{eq:common-L2-remainder} that 
		\begin{align*}
			&\left|
			\mathbb E_{\bar\mu_{N,\beta}^{(k)}}[f_{N,k}]
			-\mathbb E_{\nu_0^{\otimes k}}[f_{N,k}]
			-\frac{q_{N,2}}2
			\mathbb E_{\nu_0^{\otimes k}}
			\bigl[(\tau_k^2-k)f_{N,k}\bigr]
			\right|\\
			&\qquad\le
		\mathbb E_{\nu_0^{\otimes k}}[f_{N,k}^2]^{1/2}\E _{\nu_0^{\otimes k}} \left[\Big(L_{N,k}-1-\frac{q_{N,2}}2(\tau_k^2-k)\Big)^2\right]^{1/2}
			=\mathcal O ((kq_{N,2})^2).
		\end{align*}
		Dividing by $kq_{N,2}$ and then using $kq_{N,2}\to0$ imply the claim. 
		\end{proof}
\section{The high-temperature regime}\label{sec:high-temperature}
In this section, we consider fixed $0<\beta<1$.

\begin{lemma}\label{lem:moments}
For fixed $0<\beta<1$, the mixing law $\gamma_{N,\beta}$ converges
weakly, with convergence of all absolute moments, to
$\mathcal N(0,(1-\beta)^{-1})$. There is a constant $C=C(\beta)>0$ such
that for every $r\in\N$, we have
$\E_{\gamma_{N,\beta}}[X^{2r}]\le(Cr)^r$ and
\begin{align}\label{eq:moments31}
q_{N,2r}\le\left(\frac{Cr}{N}\right)^r\qquad \text{as well as}\qquad q_{N,2}=\frac{\beta}{(1-\beta)N}+o(N^{-1}).\end{align}
\end{lemma} 

\begin{proof}
	For each fixed $x\in\mathbb R$, Taylor
	expansion gives as $N\to\infty$
	\begin{equation}\label{eq:Taylor_FN}
		F_{N,\beta}(x)=-\frac{x^2}{2}+N\log\cosh\Big(\sqrt{\frac\beta N}x\Big)=-\frac{1-\beta}{2}x^2+O\left(\frac{x^4}{N}\right).
	\end{equation}
	The limit corresponds to a centered Gaussian density with variance $(1-\beta)^{-1}$.
Note $$0\le \log\cosh u\le u^2/2,$$ giving the upper bound $F_{N,\beta}(x)\le-\frac{1-\beta}{2}x^2$.
On the other hand, the normalizing integral in \eqref{eq:gamma} is at least
$\int_{\mathbb R}e^{-x^2/2} dx=\sqrt{2\pi}$.
It follows that
\[
\E_{\gamma_{N,\beta}}[X^{2r}]
\le\frac1{\sqrt{2\pi}}\int_{\mathbb R}
x^{2r}e^{-(1-\beta)x^2/2}dx =\frac{(2r)!}{2^r r!(1-\beta)^{r+1/2}} 
\le(Cr)^r.
\]
Dominated convergence applied to the normalizing integral and to
integrals with a fixed absolute power of $x$ proves weak convergence
and convergence of all absolute moments.

We next want to use the moment estimates for $X\sim \gamma_{N,\beta}$ to calculate the
pair correlation.  
To control the difference of $q_{N,2}$ in the form of \eqref{eq:expectation_sigmas} with $\frac{\beta}{N}\E_{\gamma_{N,\beta}}[X^2]$, we use
\begin{equation}\label{eq:tanh-quadratic-error}
|\tanh^2u-u^2|\le Cu^4,
\qquad u\in\mathbb R.
\end{equation}
Indeed, this bound follows from Taylor expansion for $|u| \le 1$ and for $|u|\ge1$ it is a consequence of $|\tanh u|\le 1$. Convergence of the second moment implies
\[
q_{N,2}=\frac\beta N\E_{\gamma_{N,\beta}}[X^2]+O\left(\frac1{N^2}\E_{\gamma_{N,\beta}}[X^4]\right)=\frac{\beta}{(1-\beta)N}+o(N^{-1}).
\]
Finally, $|\tanh u|\le|u|$ gives
$q_{N,2r}\le(\frac\beta N)^r\E_{\gamma_{N,\beta}}[X^{2r}]
\le(\frac{Cr}{N})^r$.
\end{proof}

\begin{proof}[Proof of Theorem~\ref{thm:entropy-high-temperature}]
	The case $k=1$ is trivial, hence we assume $k\ge 2$.
By Lemma~\ref{lem:moments}, pair correlation is $q_{N,2} \asymp N^{-1}$, so $k=o(N)$
implies $kq_{N,2}\to  0$. After possibly increasing the constant, the bound on
$q_{N,2r}$ yields $q_{N,2r}\le(Crq_{N,2})^r$ as in \eqref{eq:common-moment-assumption} and Proposition~\ref{prop:pair-dominance} applies. Its leading chi-square term \eqref{eq:common-chi-asymptotics} includes the count of observed pairs, i.e.
\[
\chi^2\left(\bar\mu_{N,\beta}^{(k)},\nu_0^{\otimes k}\right)
=\binom{k}{2}q_{N,2}^2+O\left((kq_{N,2})^4\right)
=\frac{\beta^2}{2(1-\beta)^2}\frac{k(k-1)}{N^2}
+o\left(\frac{k^2}{N^2}\right).
\]
The entropy is half this leading term by \eqref{eq:common-entropy-asymptotics}, proving both assertions.
\end{proof}

\begin{corollary}\label{cor:cor!}
	If $0<\beta<1$ is fixed and $k\to\infty$ with
	$k=o(N)$, then for every $g\in C_b^2(\mathbb R)$,
\begin{align}
	\mathbb E_{\bar\mu_{N,\beta}^{(k)}}
	\left[g\left(\frac{\tau_k}{\sqrt{k}}\right)\right]
	-
	\mathbb E_{\nu_0^{\otimes k}}
	\left[g\left(\frac{\tau_k}{\sqrt{k}}\right)\right]
	=
	\frac{\beta k}{2(1-\beta)N} 
	\mathbb E[g''(Z)]
	+o\left(\frac{k}{N}\right),\end{align}
	where $Z\sim\mathcal N(0,1)$.
\end{corollary} 
\begin{proof}
	Apply Corollary \ref{cor:testfunction} to
	$f_{N,k}=g(\tau_k/\sqrt{k})$ and use \eqref{eq:moments31}, we obtain 
	\begin{align}
		&\mathbb E_{\bar\mu_{N,\beta}^{(k)}}
		\left[g\left(\frac{\tau_k}{\sqrt{k}}\right)\right]
		-
		\mathbb E_{\nu_0^{\otimes k}}
		\left[g\left(\frac{\tau_k}{\sqrt{k}}\right)\right]\nonumber\\
		&\qquad=
		\frac{\beta k}{2(1-\beta)N}
		\mathbb E_{\nu_0^{\otimes k}}
		\left[\left(\frac{\tau_k^2}{k}-1\right)
		g\left(\frac{\tau_k}{\sqrt{k}}\right) 	\right]
		+o\left(\frac{k}{N}\right).\label{eq:paircor!}
	\end{align}
	Under $\nu_0^{\otimes k}$, the central limit theorem gives
	$\tau_k/\sqrt{k}\Rightarrow Z$. Since
	$\mathbb E_{\nu_0^{\otimes k}}[\tau_k^4/k^2]
	=3-2/k$, uniform integrability yields
	\[
	\mathbb E_{\nu_0^{\otimes k}}
	\left[
	g\left(\frac{\tau_k}{\sqrt{k}}\right)
	\left(\frac{\tau_k^2}{k}-1\right)
	\right]
	\longrightarrow
	\mathbb E[g(Z)(Z^2-1)]
	=\mathbb E[g''(Z)],
	\]
	where the last equality follows from Gaussian integration by (two) parts.
\end{proof}
	In particular, \eqref{eq:paircor!} reveals yet again that the first order correction of propagation of chaos is determined solely by pair correlations, through the centered squared magnetization $\tau_k^2/k-1$.
 
\section{The critical temperature}\label{sec:critical}

For $\beta=1$, we follow the same argument as in the previous section, only that the quadratic term of $\log\cosh$ cancels within $F_{N,\beta}$. The following bounds will also be used in Section \ref{sec:near-critical}.

\begin{lemma}\label{lem:h-bounds}
Let $h(u):=\log\cosh u-u^2/2$. There exist $\delta,c,C>0$ such that
\[
-Cu^4\le h(u)\le-cu^4\quad\text{for }|u|\le\delta \qquad \text{and} 
\qquad
h(u)\le-cu^2\quad\text{for }|u|\ge\delta.
\]
\end{lemma}

\begin{proof}
Taylor expansion gives $h(u)=-u^4/12+\mathcal O (u^6)$ near zero, which proves
the bounds for $u$ in a vicinity of $0$. Also, $h(0)=h'(0)=0$ and
$h''(u)=-\tanh^2u$. Thus $h(u)<0$ for $u\ne0$. Since
$h(u)/u^2\to  -1/2$ as $|u|\to  \infty$, continuity gives the stated
quadratic bound on $|u|\ge\delta$.
\end{proof}
Thus, the Taylor expansion \eqref{eq:Taylor_FN} now has leading term $-x^4/(12N)$, suggesting the rescaled mixing variable $N^{-1/4}X$ in the following critical counterpart of Lemma~\ref{lem:moments}.

\begin{lemma}\label{lem:critical-moments}
Let $X_N\sim\gamma_{N,1}$ and set $Y_N:=N^{-1/4}X_N$.
Then $Y_N$ converges weakly, with convergence of all absolute moments,
to the random variable $Y$ with density proportional to $e^{-y^4/12}$. In particular,
for every $r\in\N$, we have
$N^{r/2}q_{N,2r}\longrightarrow \E[Y^{2r}]$ and
\begin{align}\label{eq:crit_q_N_2} q_{N,2}=\frac{\alpha_{\mathrm c}}{\sqrt N}+o(N^{-1/2}).\end{align}
Moreover, there exists a constant $C>0$ such that, for all sufficiently large $N$
and all integers $1\le r< N/2$, we have $\E[Y_N^{2r}]\le(C\sqrt r)^r$
 as well as 
 \begin{align}\label{eq:crit_q_N} q_{N,2r}\le\left(\frac{C\sqrt r}{\sqrt N}\right)^r.\end{align}
\end{lemma}

\begin{proof}
The density of $Y_N$ is proportional to $e^{G_N(y)}$, where
\[
G_N(y):=F_{N,1}(N^{1/4}y)=Nh(N^{-1/4}y)=-y^4/12+\mathcal O (y^6/\sqrt N)
\]
for any fixed $y\in\R$ and, by Lemma~\ref{lem:h-bounds} 
\begin{equation}\label{eq:critical-exponent-bounds}
G_N(y)\le
\begin{cases}
-cy^4,& |y|\le\delta N^{1/4},\\
-c\sqrt N y^2,& |y|>\delta N^{1/4}.
\end{cases}
\end{equation}
In particular, $e^{G_N(y)}\le e^{-cy^4}+e^{-cy^2}$ for all $N$.
This is an integrable bound even after multiplication by any fixed
absolute power of $y$. Dominated convergence therefore applies to
the normalizing integral converging  to $\int_{\mathbb R}e^{-y^4/12}dy>0$, and to all fixed absolute-moment integrals.
It proves the asserted weak and moment convergence. 

For moments whose order grows with $N$, we use \eqref{eq:critical-exponent-bounds}. On the inner region, we have
\[
\int_{|y|\le\delta N^{1/4}}y^{2r}e^{G_N(y)}dy \le\int_{\mathbb R}y^{2r}e^{-cy^4}dy =\frac12c^{-r/2-1/4}\Gamma\left(\frac r2+\frac14\right) \le(C\sqrt r)^r,
\]
where the last inequality follows from Stirling's formula. On the outer
region, the corresponding Gaussian integral gives
\[
\int_{|y|>\delta N^{1/4}}y^{2r}e^{G_N(y)}dy
\le\int_{\mathbb R}y^{2r}e^{-c\sqrt Ny^2}dy
\le\left(\frac{Cr}{\sqrt N}\right)^r\le(C\sqrt r)^r
\]
for $r\le N$. Dividing both by the normalizing integral and increasing $C$ proves the uniform bound $\E[Y_N^{2r}]\le(C\sqrt r)^r$. Recalling the representation \eqref{eq:expectation_sigmas} for $q_{N,2r}$ and using $|\tanh u|\le|u|$, we obtain \eqref{eq:crit_q_N}. 
For fixed $r$, we 
also have convergence
\[
\left(N^{1/4}\tanh(N^{-1/4}y)\right)^{2r}
e^{G_N(y)}
\longrightarrow  y^{2r}e^{-y^4/12} ,
\] 
for all $y\in\mathbb R$. Due to
$|N^{1/4}\tanh(N^{-1/4}y)|\le|y|$ and \eqref{eq:critical-exponent-bounds}, the integrand is uniformly bounded by
$y^{2r}(e^{-cy^4}+e^{-cy^2})$.
For every fixed $r$, it follows again from dominated convergence
\[
N^{r/2}q_{N,2r} \longrightarrow 
\frac{\displaystyle\int_{\mathbb R}y^{2r}e^{-y^4/12} dy}
{\displaystyle\int_{\mathbb R}e^{-y^4/12} dy}.
\]
Taking $r=1$ yields \eqref{eq:crit_q_N_2}. 
\end{proof}

The restriction $r\le N/2$ will be sufficient, since it covers all coefficients needed
in the expansion \eqref{eq:LNWalsh}. 

\begin{proof}[Proof of Theorem~\ref{thm:critical-entropy}]
	The case $k=1$ is trivial, hence we assume $k\ge 2$. Lemma~\ref{lem:critical-moments} gives
$q_{N,2}\asymp N^{-1/2}$, so $k=o(\sqrt N)$ implies
$kq_{N,2}\to  0$. Since $\sqrt r\le r$ for $r\ge1$, the uniform
bound \eqref{eq:crit_q_N} implies
$q_{N,2r}\le(Crq_{N,2})^r$.
Proposition~\ref{prop:pair-dominance} therefore applies. Its pair
term satisfies
\[
\binom{k}{2}q_{N,2}^2
=\frac{\alpha_{\mathrm c}^2}{2}\frac{k(k-1)}{N}
+o\left(\frac{k^2}{N}\right),
\]
and the higher-degree error is $\mathcal O (k^4/N^2)=o(k^2/N)$.
The entropy is half the leading chi-square term.
\end{proof}

\section{Approaching the critical temperature}\label{sec:near-critical}

We now consider $\beta_N\uparrow1$ and set $\lambda_N=\sqrt N(1-\beta_N)$. 
We first collect properties of $\alpha(\lambda)$ and uniform estimates on the spin correlations.

\begin{lemma}\label{lem:alpha-properties}
The function $\alpha$ is continuous and strictly decreasing, with
$\alpha(0)=\alpha_{\mathrm c}$ and
$\alpha(\lambda)\sim\lambda^{-1}$ as $\lambda\to  \infty$.
There are constants $c,C>0$ such that
\begin{equation}\label{eq:alpha-comparison}
\frac{c}{1+\lambda}\le\alpha(\lambda)\le\frac{C}{1+\lambda},
\qquad \text{for} \quad\lambda\ge0.
\end{equation}
\end{lemma}

\begin{proof}
Let $Y_\lambda$ have density proportional to
$\exp(-\lambda y^2/2-y^4/12)$. By the definition of $\alpha$ in
Theorem~\ref{thm:near-critical}, we have
$\alpha(\lambda)=\E[Y_\lambda^2]$ and differentiation under the
integral gives
\[
\alpha'(\lambda)=-\frac12\left(\E[Y_\lambda^4]-\E[Y_\lambda^2]^2\right)
=-\frac12\Var(Y_\lambda^2)<0,
\]
where the strict inequality follows because $Y_\lambda^2$ is not
constant. Continuity follows from dominated convergence, and
$\alpha(0)=\alpha_{\mathrm c}$ by definition. The representation in
terms of $\Gamma$-functions in \eqref{eq:alphac} follows by substituting $t=y^4/12$ in the
normalizing and second-moment integrals. For $\lambda>0$, the change of variables $z=\sqrt\lambda y$
gives
\[
\lambda\alpha(\lambda)
=\frac{\displaystyle\int_{\mathbb R}z^2e^{-z^2/2-z^4/(12\lambda^2)} dz}
{\displaystyle\int_{\mathbb R}e^{-z^2/2-z^4/(12\lambda^2)} dz}
\longrightarrow  1 \qquad\text{as }\lambda\to  \infty
\]
by dominated convergence.
Positivity and continuity on bounded intervals now give \eqref{eq:alpha-comparison}. 
\end{proof}

The following lemma provides the correlation estimates
needed to apply Proposition~\ref{prop:pair-dominance}.

\begin{lemma}\label{lem:uniform-near-critical-moments}
Let $\beta_N\uparrow1$ and $\lambda_N=\sqrt N(1-\beta_N)$. Then
\begin{equation}\label{eq:qN2_betaN}
q_{N,2}=\frac{\alpha(\lambda_N)}{\sqrt N}(1+o(1)).
\end{equation}
Moreover, there is a constant $C>0$ such that, for all sufficiently
large $N$ and all integers $1\le r\le N/2$,
\begin{equation}\label{eq:bound_q2r_betaN}
q_{N,2r}\le N^{-r/2}
\left(C\min\left\{\sqrt r,\frac{r}{1+\lambda_N}\right\}\right)^r.
\end{equation}
\end{lemma}

\begin{proof}
Let $X_N\sim\gamma_{N,\beta_N}$ and $Y_N=N^{-1/4}X_N$.
Its density is proportional to $e^{G_N(y)}$, with
\begin{equation}\label{eq:near-critical-exponent}
	G_N(y):=F_{N,\beta_N}(N^{1/4}y)
	=-\frac{\lambda_N}{2}y^2
	+Nh(\sqrt{\beta_N}N^{-1/4}y).
\end{equation}
We may assume $\beta_N\ge1/2$ by taking $N$ sufficiently large.
	
	On $|y|\le(1+\lambda_N)^{-1/2}$,
	Lemma~\ref{lem:h-bounds} directly gives $G_N(y)\ge-C$, hence
	\begin{equation}\label{eq:near-critical-normalization}
		\int_{\mathbb R}e^{G_N(y)}\,dy
		\ge\frac{c}{\sqrt{1+\lambda_N}}.
	\end{equation}
	
	We shall now establish the following uniform moment bound 
\begin{equation}\label{eq:combined-moment-bound}
	\E[Y_N^{2r}]=\frac{\int_{\mathbb R}|y|^{2r}e^{G_N(y)} dy}
	{\int_{\mathbb R}e^{G_N(y)} dy}
	\le\left(C\min\left\{\sqrt r,
	\frac{r}{1+\lambda_N}\right\}\right)^r,
	\qquad 1\le r\le N/2.
\end{equation}
	To this end, let $\widetilde Y_N$ have density proportional to $	\exp\big(Nh(\sqrt{\beta_N}N^{-1/4}y)\big)$.
	Since $\beta_N\ge 1/2$, the estimates from
	Lemma~\ref{lem:h-bounds} give the same moment bounds as in Lemma~\ref{lem:critical-moments}, that is
$$\E[\widetilde Y_N^{2r}]
	\le\beta_N^{-r}(C\sqrt r)^r
	\le(C\sqrt r)^r.$$
	The law of $Y_N$ is obtained by tilting the density of
	$\widetilde Y_N$ with $e^{-\lambda_Ny^2/2}$.
	Since $|y|^{2r}$ is increasing and
	$e^{-\lambda_Ny^2/2}$ is decreasing in $|y|$,
	the covariance of $|\widetilde Y_N|^{2r}$ and $e^{-\lambda_N\widetilde Y_N^2/2}$ 
	is nonpositive. Hence, we obtain the first bound of \eqref{eq:combined-moment-bound},
\begin{align}\label{eq:firstbound}
	\E[Y_N^{2r}]
	=\frac{\E[\widetilde Y_N^{2r}
		e^{-\lambda_N\widetilde Y_N^2/2}]}
	{\E[e^{-\lambda_N\widetilde Y_N^2/2}]}
	\le\E[\widetilde Y_N^{2r}]
	\le(C\sqrt r)^r. \end{align}

	If $\lambda_N\ge1$, then $h\le0$ and \eqref{eq:near-critical-normalization} yields
	\begin{align*}
		\E[Y_N^{2r}]
		\le
		\frac{\int_{\mathbb R}
			y^{2r}e^{-\lambda_Ny^2/2}dy}
		{c(1+\lambda_N)^{-1/2}}
		\le C\sqrt{\frac{1+\lambda_N}{\lambda_N}}
		\left(\frac{Cr}{\lambda_N}\right)^r
		\le\left(\frac{Cr}{1+\lambda_N}\right)^r.
	\end{align*} 
	For $\lambda_N<1$, the same bound follows from \eqref{eq:firstbound},
	since $\sqrt r\le 2r/(1+\lambda_N)$.
	Combining the two estimates, we obtain our claim \eqref{eq:combined-moment-bound}.
	
	Finally, the representation \eqref{eq:expectation_sigmas} of $q_N$
	and $|\tanh u|\le|u|$ imply
	\[
	q_{N,2r}=\E\Big[\tanh^{2r}(\sqrt{\beta_N}N^{-1/4}Y_N)\Big]
	\le\beta_N^rN^{-r/2}\E[Y_N^{2r}]
	\]
and	\eqref{eq:combined-moment-bound}  proves \eqref{eq:bound_q2r_betaN}.
 
Now we turn to analyzing $q_{N,2}$.
To determine its asymptotic behavior, we first
prove the relative comparison
\begin{equation}\label{eq:second-moment-comparison}
	\E[Y_N^2]=\alpha(\lambda_N)(1+o(1)).
\end{equation}
We  use the subsequence criterion for convergence: it suffices
to show that every subsequence of $\E[Y_N^2]/\alpha(\lambda_N)$ has a further subsequence
converging to $1$.

Starting from an arbitrary subsequence, we select a further
subsequence along which $\lambda_N\to  \lambda\le\infty$. In the two cases below, we prove
\eqref{eq:second-moment-comparison} along the selected
subsequence.

Assume first $\lambda_N\to  \lambda<\infty$, then for each fixed $y$,
$G_N(y)\to   -\frac\lambda2y^2-\frac1{12}y^4$.
Using the bounds from Lemma~\ref{lem:h-bounds} in the
same way as in the proof of
Lemma~\ref{lem:critical-moments}, we obtain an integrable
majorant independent of $N$ for both $e^{G_N(y)}$ and
$y^2e^{G_N(y)}$. Dominated convergence applied to the
normalizing and second-moment integrals gives
$\E[Y_N^2]\to  \alpha(\lambda)$. Since
$\alpha(\lambda_N)\to  \alpha(\lambda)>0$, we obtain 
\eqref{eq:second-moment-comparison}.

If $\lambda_N\to  \infty$, we need to rescale once more and set
$Z_N:=\sqrt{\lambda_N} Y_N$. This keeps the second moment on a non-degenerate scale. The density of $Z_N$ has exponent
\[
G_N\Big(\frac{z}{\sqrt{\lambda_N}}\Big)
=-\frac{z^2}{2}
+Nh\Big(\frac{\sqrt{\beta_N} z}{N^{1/4}\sqrt{\lambda_N}}\Big).
\]
For fixed $z$, the second term is
\[
-\frac{\beta_N^2z^4}{12\lambda_N^2}
+O\Big(\frac{z^6}{\sqrt N \lambda_N^3}\Big) \longrightarrow  0.
\]
Now, in a similar fashion as before, we write the second moment of $Z_N$ as
\begin{align}\label{eq:numden}
\E[Z_N^2]=
\frac{\displaystyle\int_{\mathbb R}z^2e^{G_N(z/\sqrt{\lambda_N})} dz}
{\displaystyle\int_{\mathbb R}e^{G_N(z/\sqrt{\lambda_N})} dz}.
\end{align}
The pointwise convergence established above gives
$e^{G_N(z/\sqrt{\lambda_N})}\to   e^{-z^2/2}$.
Moreover, $h\le0$ implies the bounds
\[
0\le e^{G_N(z/\sqrt{\lambda_N})}\le e^{-z^2/2}
\qquad  \text{and} \qquad 
0\le z^2e^{G_N(z/\sqrt{\lambda_N})}\le z^2e^{-z^2/2}.
\]
Both functions on the right are integrable. We can therefore apply dominated convergence separately to the denominator
and numerator of \eqref{eq:numden}, each converging to $\sqrt{2\pi}$.
Hence, the quotient tends to one. Thus
$
\lambda_N\E[Y_N^2]=\E[Z_N^2]\longrightarrow 1.
$
Lemma~\ref{lem:alpha-properties} also gives
$\lambda_N\alpha(\lambda_N)\to   1$ and therefore
\[
\frac{\E[Y_N^2]}{\alpha(\lambda_N)} =
\frac{\lambda_N\E[Y_N^2]} {\lambda_N\alpha(\lambda_N)}
\longrightarrow  1.
\]
This proves \eqref{eq:second-moment-comparison} in the case
$\lambda_N\to  \infty$, and the claim \eqref{eq:second-moment-comparison} follows.


Let us apply this comparison to determine $q_{N,2}$.
Replacing $\tanh^2(u)$ by $u^2$, we make the same error as in \eqref{eq:tanh-quadratic-error} with
$u=\sqrt{\beta_N}N^{-1/4}Y_N$, then take expectations, and multiply by $\sqrt N$, yielding
\[
\left|\sqrt N q_{N,2}-\beta_N\E[Y_N^2]\right|
\le \frac{C\beta_N^2}{\sqrt N}\E[Y_N^4]
\le \frac{C}{\sqrt N}\E[Y_N^4],
\]
where the last inequality uses $\beta_N\le1$.

To compare this error with $\alpha(\lambda_N)$, we use
\eqref{eq:combined-moment-bound} with $r=2$ and the lower
bound in \eqref{eq:alpha-comparison}. After increasing $C$
if necessary, these give
$
\E[Y_N^4]\le\frac{C}{(1+\lambda_N)^2}
$
and $\alpha(\lambda_N)\ge\frac{c}{1+\lambda_N}$.
Dividing the error estimate by $\alpha(\lambda_N)$
therefore gives
\[
\left|\frac{\sqrt N q_{N,2}}{\alpha(\lambda_N)}-\beta_N\frac{\E[Y_N^2]}{\alpha(\lambda_N)}\right|
\le\frac{C \E[Y_N^4]}{\sqrt N \alpha(\lambda_N)}\le
\frac{C^2}{c\sqrt N(1+\lambda_N)} \longrightarrow  0.
\]
Finally, $\beta_N\to   1$ and \eqref{eq:second-moment-comparison} imply
$\beta_N\frac{\E[Y_N^2]}{\alpha(\lambda_N)} \longrightarrow  1$.
Hence $\sqrt N q_{N,2}/\alpha(\lambda_N)\to   1$,  which is \eqref{eq:qN2_betaN}.
\end{proof}

The estimates now allow to verify the conditions of
Proposition~\ref{prop:pair-dominance}. 

\begin{proof}[Proof of Theorem~\ref{thm:near-critical}]
The stated properties of $\alpha$, including its range, were proved
in Lemma~\ref{lem:alpha-properties}. 	For $k=1$, both divergences vanish identically.
For $k\ge2$, Lemma~\ref{lem:uniform-near-critical-moments} gives
$
q_{N,2}=\frac{\alpha(\lambda_N)}{\sqrt N}(1+o(1)),
$
so the assumed block-size condition implies $kq_{N,2}\to   0$.
On the other hand, from Lemma~\ref{lem:uniform-near-critical-moments} 
we obtain
\[
q_{N,2r}\le N^{-r/2}\left(C\min\left\{\sqrt r,\frac{r}{1+\lambda_N}\right\}\right)^r\le
\left(\frac{Cr}{\sqrt N(1+\lambda_N)}\right)^r.
\]
Moreover, the asymptotic formula for $q_{N,2}$
and the lower bound in \eqref{eq:alpha-comparison} imply,
for all sufficiently large $N$,
$q_{N,2}\ge\frac{\alpha(\lambda_N)}{2\sqrt N}\ge\frac{c}{2\sqrt N(1+\lambda_N)}$ and therefore
\[
q_{N,2r}\le \left(\frac{2C}{c} r q_{N,2}\right)^r.
\]
Since $2/c$ is independent of $N$ and $r$, we can absorb
it into $C$. This verifies
\eqref{eq:common-moment-assumption}
for all sufficiently large $N$ and we can apply Proposition~\ref{prop:pair-dominance}. Its leading chi-square term is
\[
\binom{k}{2}q_{N,2}^2
=\frac{\alpha(\lambda_N)^2}{2}\frac{k(k-1)}{N}
+o\left(\frac{k^2\alpha(\lambda_N)^2}{N}\right).
\]
The higher-degree error is
$\mathcal O ((k\alpha(\lambda_N)/\sqrt N)^4)$, which is negligible compared
with $k^2\alpha(\lambda_N)^2/N$ by assumption. This proves
\eqref{eq:chi2window}, and
\eqref{eq:common-entropy-asymptotics} proves \eqref{eq:KLwindow}.
\end{proof}

%

\end{document}